\documentclass[11pt]{article}

\usepackage{amssymb,latexsym}
\usepackage{epsfig}
\usepackage{eufrak}
\usepackage{amsmath}
\usepackage{mathrsfs}
\usepackage{color}

\newtheorem{theorem}{Theorem}
\newtheorem{lemma}{Lemma}

\newtheorem{remark}{Remark}

\newcommand{\be}{\begin{equation}}
\newcommand{\ee}{\end{equation}}
\newcommand{\bee}{\begin{eqnarray*}}
\newcommand{\eee}{\end{eqnarray*}}
\newcommand{\bel}{\begin{eqnarray}}
\newcommand{\eel}{\end{eqnarray}}
\newcommand{\bec}{\begin{cases}}
\newcommand{\eec}{\end{cases}}
\newcommand{\bem}{\begin{bmatrix}}
\newcommand{\eem}{\end{bmatrix}}

\newcommand{\la}{\label}
\newcommand{\li}{\left}
\newcommand{\ri}{\right}

\newcommand{\DEF}{\stackrel{\mathrm{def}}{=}}

\newcommand{\ep}{\epsilon}

\newcommand{\lm}{\lambda}

\newcommand{\si}{\sigma}

\newcommand{\de}{\delta}

\newcommand{\se}{\theta}
\newcommand{\Se}{\Theta}

\newcommand{\ze}{\zeta}
\newcommand{\al}{\alpha}
\newcommand{\ba}{\beta}

\newcommand{\ro}{\rho}
\newcommand{\ka}{\kappa}

\newcommand{\f}{\frac}
\newcommand{\sq}{\sqrt}
\newcommand{\cd}{\cdots}

\newcommand{\qu}{\quad}
\newcommand{\qqu}{\qquad}
\newcommand{\fa}{\forall}

\newcommand{\bb}{\mathbb}

\newcommand{\wh}{\widehat}

\newcommand{\mrm}{\mathrm}

\newcommand{\tx}{\text}

\newcommand{\iy}{\infty}

\newcommand{\bed}{\begin{description}}
\newcommand{\eed}{\end{description}}
\newcommand{\bei}{\begin{itemize}}
\newcommand{\eei}{\end{itemize}}
\newcommand{\ben}{\begin{enumerate}}
\newcommand{\een}{\end{enumerate}}
\newcommand{\bib}{\bibitem}
\newcommand{\beL}{\begin{lemma}}
\newcommand{\eeL}{\end{lemma}}
\newcommand{\beT}{\begin{theorem}}
\newcommand{\eeT}{\end{theorem}}
\newcommand{\sect}{\section}

\newcommand{\bpf}{\begin{pf}}
\newcommand{\epf}{\end{pf}}
\newcommand{\bsk}{\bigskip}

\newcommand{\pfbox}{\hfill\mbox{$\Box$}}

\newenvironment{pf}{\paragraph*{Proof{\rm.}}}{\pfbox\bigskip}

\begin{document}

\title{{\bf Risk Analysis in Robust Control --- Making the Case for
Probabilistic Robust Control}
\thanks{This research is supported in part by the US Air Force.
\newline The authors are with
Department of Electrical and Computer Engineering, Louisiana State
University, Baton Rouge, LA 70803; Email: \{chan,kemin,
aravena\}@ece.lsu.edu, Tel: (225)578-\{8961, 5533,5537\}, and Fax:
(225) 578-5200. }}

\author{Xinjia Chen, Jorge L. Aravena and Kemin Zhou}

\date{June 2007}

\maketitle

\begin{abstract}

This paper offers a critical view of the ``worst-case"  approach
that is the cornerstone of robust control design. It is our
contention that a blind acceptance of worst-case scenarios may
lead to designs that are actually more dangerous than designs
based on probabilistic techniques with a built-in risk factor. The
real issue is one of modeling. If one accepts that no
mathematical model of uncertainties is perfect then a
probabilistic approach can lead to more reliable control even if
it cannot guarantee stability for all possible cases. Our
presentation is based on case analysis. We first establish that
worst-case is not necessarily ``all-encompassing." In fact, we
show that for some uncertain control problems to have a
conventional robust control solution it is necessary to make
assumptions that leave out some feasible cases. Once we establish
that point, we argue that it is not uncommon for the risk of
unaccounted cases in worst-case design to be greater than that of
the accepted risk in a probabilistic approach. With an
example, we quantify the risks and show that worst-case can be
significantly more risky.  Finally, we join our analysis with existing results
on computational complexity and probabilistic robustness to argue that the deterministic
worst-case analysis is not necessarily the better tool.
\end{abstract}

\bsk

\section{Introduction}

In recent years, a number of researchers have proposed probabilistic control
methods as alternatives for overcoming the computational
complexity and conservatism of deterministic
worst-case robust control framework (e.g.,
\cite{bai}--\cite{SB}, \cite{kov}--\cite{Wang2} and the references
therein). The philosophy of probabilistic control theory is to
sacrifice extreme cases of uncertainty. Such paradigm has lead to
the novel concepts of probabilistic robustness margin and
confidence degradation function (e.g., \cite{BLT}). Despite
the claimed advantages of probabilistic approach, the
deterministic worst-case approach remains dominating for design and analysis
purposes. It
is a common contention that a probabilistic design is more risky
than a worst-case design. Such a contention may have been the main
obstacle preventing the wide acceptance of the probabilistic
paradigm, especially in the development of highly reliable
systems. When
referring to the probabilistic approach, a cautious warning is
usually attached. Statements like ``if one is
willing to accept a small risk of performance violation'' can be found in a
number of robust control papers. A
typical argument is that the worst-case method takes every case of
``uncertainty'' into account and is certainly the most safe,
while the probabilistic method considers only most of the
instances of ``uncertainty'' and, hence, is more risky.

We illustrate
with two very simple cases that a worst case design may not necessarily
consider all possible cases. Our purpose here is to make the point that
the basic issue is one of modeling and as such it is
never perfect. Worst-case scenario may mean ``the worst case that
we can imagine," or ``the
worst-case that we can afford to consider to have a robust
solution."

In practice, the coefficients of a linear model are complex
functions of physical parameters.  Even if the physical parameters
are bounded in a narrow interval, the variations of the
coefficients can be fairly large.  A simple example is provided by the
process of discharge of a cylindrical tank. The basic nonlinear model is of the
form
\[
A\frac{dH}{dt}+\rho\sqrt{H}=Q_i
\]
where $A$ is the tank cross section, $H$ the height of liquid inside the
tank, $Q_i$ the
volumetric input flow rate
and $\rho$ the hydraulic resistance in the discharge. Linearizing in the
neighborhood of a steady
state operating point, ($\overline{H},\overline{Q}_i$), satisfying $\rho
\sqrt{\overline{H}}
=\overline{Q_i}$ one
obtains the linear model
\[
\frac{dh}{dt}+\frac{\rho}{2A\sqrt{\overline{H}}}h = \frac{q_i}{A}
\]
with $h=H-\overline{H}, q_i=Q_i-\overline{Q}_i$. Clearly, the parameter $a=
\frac{\rho}{2A\sqrt{\overline{H}}}$ takes values
in the interval
($0,~\infty$). Hence, any design assuming bounded uncertainties for the
parameter, $a$, cannot
include all possible heights for the liquid.

For the second case consider a first order system of the form
\[
G(s) = \f{q}{s - p}
\]
with uncertain  parameters $q$ and $p$. Assuming a unity feedback and
controller of the form
\[
C(s)=\frac{K}{s+a}, \qu K > 0, \qu a > 0
\]
the closed loop system becomes
\[
T(s)=\frac{K q}{s^2+(a-p)s+K q - a p}.
\]
The controller will robustly stabilize the plant if
\[
p  <  a, \qqu q  >  \frac{a p}{K}.
\]
It is clear that for any finite controller gain $K$, there
exist a range of
values of $q$ where the closed-loop system
is unstable. The designer of a worst-case controller would be faced with the
choice of selecting
a different controller
structure or assuming, based on other considerations, that a neighborhood
of $q = 0$ can be
excluded from the design. With the next result we develop this point in a
more
general form.

\subsection{Uncertainties in Modeling Uncertainties}

In many practical situations of worst-case design one models
uncertainties as bounded random variables. The issue of selecting
the bounds is not trivial and is, oftentimes, not addressed in detail. The
following theorem
shows that, regardless of the assumed size of the uncertainty
set, a worst-case robust controller actually can always fail.
Hence, if there are ``uncertainties in modeling the uncertainties"
it may be better to model them as random variables varying from
$-\infty$ to $\infty$ in order to pursue ``worst-case'' in a
strict sense.

\beT Let the transfer function of the uncertain plant be
\[
G(s) = \f{ \sum_{i=0}^\ell  \ba_i s^{\ell- i}  } { \sum_{i=0}^\ka
\al_i s^{\ka- i}  }, \qqu \al_0 = 1, \qqu \ell \leq \ka.
\]
Assume that for a given finite uncertainty range in the parameters
$\al_i, \; i = 1 , \cd, \ka$ and $\ba_j, \; j = 0, 1, \cd, \ell$,
there exists a controller of the form
\[
C(s) = \f{ \sum_{i=0}^m  b_i s^{m- i}  } { \sum_{i=0}^n a_i
s^{n- i}  },  \qqu a_0 = 1, \qqu b_0 \neq 0, \qqu m \leq n
\]
which robustly stabilizes the system.   Then, there always exists
a value of parameter $\al_i$ or $\ba_j$, outside the assumed
uncertainty range and which will make the closed-loop system
unstable. \eeT

\bpf The characteristic equation of the closed-loop system is
\[
\li ( \sum_{i=0}^n a_i s^{n- i} \ri ) \; \li ( \sum_{j=0}^\ka
\al_j s^{\ka - j}  \ri ) + \li ( \sum_{i=0}^m  b_i s^{m- i}  \ri )
\; \li ( \sum_{j=0}^\ell  \ba_j s^{\ell- j}  \ri ) = 0,
\]
which can be written as
\[
\sum_{\tau = 0}^{n + \ka}  \sum_{i +j = \tau \atop{ 0 \leq i \leq
n  \atop{0 \leq j \leq \ka  } } }  (a_i \; \al_j) \; s^{n+ \ka -
\tau}  + \sum_{\iota = 0}^{m + \ell}  \sum_{i +j = \iota \atop{ 0
\leq i \leq m  \atop{0 \leq j \leq \ell   } } }  (b_i \; \ba_j) \;
s^{m + \ell - \iota} = 0.
\]
For $1 \leq \tau \leq \ka$,  the coefficient of $s^{n+ \ka -
\tau}$ is
\[
\sum_{i +j = \tau \atop{ 0 \leq i \leq n  \atop{0 \leq j \leq \ka
} } }  (a_i \; \al_j)  +   \sum_{i +j = \tau + m + \ell - (n +
\ka) \atop{ 0 \leq i \leq m \atop{0 \leq j \leq \ell   } } }  (b_i
\; \ba_j) = \al_\tau + \xi
\]
where \bee \xi & = & \sum_{i +j = \tau \atop{ 1 \leq i \leq n
\atop{0 \leq j \leq \ka } } }  (a_i \; \al_j)  +   \sum_{i +j =
\tau + m + \ell - (n + \ka) \atop{ 0 \leq i \leq m \atop{0 \leq j
\leq \ell } } } (b_i \; \ba_j)\\
& = & \sum_{j = \max (0, \tau - n)  }^{\tau -1}  (a_{\tau - j} \;
\al_j) + \sum_{i +j = \tau + m + \ell - (n + \ka) \atop{ 0 \leq i
\leq m \atop{0 \leq j \leq \ell } } } (b_i \; \ba_j) \eee is
independent of $\al_\tau$ and depends on the other uncertainties.
It follows that the system will be unstable if
\be
\la{con3}
\al_\tau \leq - \xi.
\ee

In a similar manner, for $0 \leq \iota \leq \ell$,  the
coefficient of $s^{m+ \ell - \iota}$ is
\[
\sum_{i +j = \iota \atop{ 0 \leq i \leq m  \atop{0 \leq j \leq
\ell   } } }  (b_i \; \ba_j)  + \sum_{i +j = \iota + n + \ka - (m
+ \ell) \atop{ 0 \leq i \leq n  \atop{0 \leq j \leq \ka  } } }
(a_i \; \al_j) = b_0 \; \ba_\iota + \ze
\]
where \bee \ze & = & \sum_{i +j = \iota \atop{ 1 \leq i \leq m
\atop{0 \leq j \leq \ell   } } }  (b_i \; \ba_j)  + \sum_{i +j =
\iota + n + \ka - (m + \ell) \atop{ 0 \leq i \leq n  \atop{0 \leq
j \leq \ka  } } } (a_i \; \al_j)\\
& = & \sum_{ j = \max(0, \iota - m) }^{ \iota - 1 } (b_{\iota - j}
\; \ba_j) + \sum_{i +j = \iota + n + \ka - (m + \ell) \atop{ 0
\leq i \leq n  \atop{0 \leq j \leq \ka  } } } (a_i \; \al_j) \eee
is independent of $\ba_\iota$ and depends on the other
uncertainties. It follows that the system will be unstable if
\be
\la{con8}
b_0 \; \ba_\iota + \ze \leq 0.
\ee
In the case that $b_0 > 0$, the system is unstable
when $\ba_\iota \leq - \f{\ze}{b_0}$.  In the case that $b_0 < 0$,
the system is unstable when $\ba_\iota \geq - \f{\ze}{b_0}$.

Let $\cal{B}$ be the uncertainty bounding set {\em assumed} by the designer.
Let $\cal{D}$ be the set of all values of uncertainties for which
the controller stabilizes the system. Obviously, $\cal{B} \subseteq
\cal{D}$.
From the stability conditions (\ref{con3}) and (\ref{con8}),
we can conclude that $\cal{D}$ must be bounded.
Taking into account the stability conditions and the fact that $\cal{D}$ is
bounded,
we can see that there exist values of
the parameters $\al_i, \; 1 \leq i \leq \ka$ or $\ba_j, \; 0 \leq
j \leq \ell$ that fall outside the domain $\cal{D}$
(of course they fall outside the assumed uncertainty range $\cal{B}$) and
make
the
system unstable.

\epf

\begin{remark}
A proof for an equivalent result for a multi-variable plant may be feasible.
However, the
following general argument conveys the idea about the limitations in worst
case design.  Assume that each instance of uncertainty is an element of $n$-dimensional vector space $E^n$.
Let $G(s, q), \; q \in\cal{B}\subset E^n$ be the model for an uncertain
plant.
Assume that there exists a
controller $C_w$ that satisfies the robustness requirements for all $q\in
\cal{B}$. Define now
as $\cal{D}$ the set of all values of the parameter $q$ where the controller
$C_w$ satisfies the
robustness requirements. Clearly $\cal{D}_w\supset\cal{B}$  but unless
$\cal{D}_w=E^n$ there always exist values of the parameter $q$ where the
controller does not
perform. The worst-case design ignores these cases as impossible. Our
contention is that the
modeling of uncertainties (the set $\cal{B}$) may not include all cases that
could occur and it may be
better to accept a risk from the onset of the design.
\end{remark}

\section{Accepting Risk Can Be Less Risky}
The previous result makes, very strongly, the
point that worst-case modeling is not ``all-encompassing" and
therefore it has some risk associated to it. In this section we offer first
a more
formal description of the problem and argue that a
probabilistic approach may easily lead to more reliable designs. The next
section
uses a case study to quantify the actual risks of both approaches.
\subsection{Designing with Uncertain Uncertainties}
We incorporate the fact that modeling is never exact by
postulating an uncertainty set, $\cal{U}$ and a bounding set,
$\cal{B}$, that models the uncertainties. The actual relationship
between the two sets is not known. The worst-case design finds a
controller $C_w$ to guarantee every uncertainty instance $q \in
\cal{B}$. The probabilistic design seeks a controller $C_p$ to
guarantee most uncertainty instances $q \in \cal{B}$. Formally we
can define the following relevant subsets \bee
\cal{M}&=& \cal{U}\cap \cal{B},\\
\cal{N}&=& \overline{\cal{U}}\cap \cal{B},\\
\cal{E}&=& \cal{U}\cap \overline{\cal{B}}.\eee
Here $\overline{X}$ denotes the complementary set of $X$.
Clearly, $\cal{M}$ contains those uncertainties that are modeled
while the set $\cal{N}$ describes modeled uncertainties that
never occur and $\cal{E}$ describes the unmodeled uncertainties.
The existence of these last two sets creates either inefficiencies
or risks in the worst case design. To see this, consider the
extreme situation where the designed robust controller guarantees
the robustness requirement only for instances in $\cal{B}$. The
controller is over-designed because it deals with situations that
cannot occur and it has the risk of failure if an instance in the
set $\cal{E}$ arises.

Having established the fact that a
robust control design can be risky, we now argue that
probabilistic design can actually be less risky. As an added
benefit, it is known that many worst-case synthesis problems are either not
tractable, or have known solutions which are unduly
conservative and implementation expensive. But when using a
probabilistic method, the previously intractable problems may
become solvable, the conservatism may be overcame, and high
performance controller with simple structure may be obtained.

For brevity, we use notation $C^V$  to represent the statement
that ``the robustness requirement is violated for the system with
controller $C$''.  The subindex $w$ will refer to worst-case
design while $p$ will refer to probabilistic design. Note that
the risk of a probabilistic design is
\[
P_e^p = \Pr \{ C^V_p \mid q \in \cal{M}  \} \; \Pr\{ q \in
\cal{M}  \} + \Pr \{ C^V_p  \mid q \in \cal{E}  \} \; \Pr \{ q
\in \cal{E}  \}.
\]
While the risk of a worst-case design is
\[
P_e^w = \Pr \{ C^V_w  \mid q \in \cal{E}  \} \; \Pr \{ q \in
\cal{E}  \}.
\]
Hence the ratio of risks will be \be
\la{eeq}
\f{ P_e^p } { P_e^w }  =
\frac{\Pr\{C^V_p \mid q \in \cal{E}\}}{\Pr\{C^V_w \mid q \in
\cal{E}\}}+ \f{ \Pr \{ C^V_p
\mid q \in \cal{M}  \}   \Pr \{ q \in \cal{M}  \} }{ \Pr \{ C^V_w  \mid q
\in \cal{E} \} \Pr\{q\in \cal{E}\}  }. \ee
The first term is related to the
performance of both types of controllers outside the design
region. The behavior of the worst-case design in this region is
of no concern to the designer, after all it ``never gets there.''
All the design effort is placed in assuring performance over the
set $\cal{B}$. Hence we can reasonably expect $\Pr \{ C^V_w  \mid q
\in \cal{E} \}$ to be high. In fact, if the set $\cal{N}$, of impossible
situations
included  in the design, is large then $\Pr \{ C^V_w  \mid q
\in \cal{E} \}$ could be close to one and the first term in the right-hand
side of (\ref{eeq})
can easily be less than some number $\lm \in (0,1)$. The second
term contains the factor $\Pr\{C^V_p \mid q \in \cal{M}\}$ which
is under the probabilistic designer and is a measure of the
accepted risk. It is reasonable to expect that this risk is less than
$\Pr\{q\in \cal{E}\}$
so that \bee \frac{\Pr\{C^V_p
\mid q \in \cal{M}\}\Pr\{q\in \cal{M}\}}{\Pr\{q\in \cal{E}\}}< (1 - \lm) \;
\Pr \{ C^V_w  \mid q
\in \cal{E} \}.
\eee
and consequently $\f{ P_e^p } { P_e^w }  < 1$.  Factoring in a poor
performance for the worst-case design outside
of the modeled region one can see that the risk of the
probabilistic design can become smaller than the risk of
a worst-case design with unmodeled uncertainties.

From a different point of view, many experiments of performance
degradation of probabilistic designs indicate a fairly flat
characteristic. If the unmodeled uncertainty set $\cal{E}$ is
relatively small then one could argue that \bee \Pr\{C^V_p \mid
q\in\cal{E}\} \approx \Pr\{C^V_p \mid q\in\cal{M}\} \eee and the
ratio of risks is approximately given by \be \f{ P_e^p } { P_e^w
}  \approx  \frac{ \Pr\{C^V_p \mid q\in\cal{M} \}}{\Pr \{ C^V_w
\mid q \in \cal{E} \}\Pr\{q\in \cal{E}\}  }. \ee The numerator is
under the control of the designer in a probabilistic approach
while the denominator has not even been considered  as existing
in a worst-case design.

\sect{Comparing Worst-Case and Probabilistic Designs}
The problem of quantifying the differences in performance between a
worst-case design and a
probabilistic design is extremely difficult in general.
In this section we use a case study to quantify the risks and
demonstrate  that it is {\it not uncommon} for a
probabilistic
controller to be (highly) more reliable than a worst-case controller. We
postulate that if the
result holds for simple systems
then it is also likely for more complex situations.

\bsk

Consider a feedback system as follows.
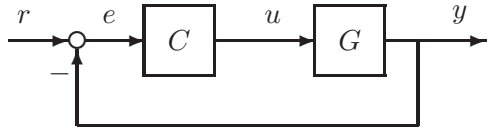
\begin{figure}[htb]
\centering{
\setlength{\unitlength}{0.0006in}%
\begin{picture}(4224,1074)(1189,-3523)
\thicklines \put(1801,-2761){\circle{150}}
\put(2401,-3061){\framebox(600,600){$C$}}
\put(1201,-2761){\vector( 1, 0){525}} \put(1876,-2761){\vector( 1,
0){525}} \put(3001,-2761){\vector( 1, 0){900}}
\put(3901,-3061){\framebox(600,600){$G$}}
\put(4501,-2761){\vector( 1, 0){900}} \put(4801,-2761){\line(
0,-1){750}} \put(4801,-3511){\line(-1, 0){3000}}
\put(1801,-3511){\vector( 0, 1){675}}
\put(1276,-2611){\makebox(0,0)[lb]{$r$}}
\put(2026,-2611){\makebox(0,0)[lb]{$e$}}
\put(5101,-2611){\makebox(0,0)[lb]{$y$}}
\put(3451,-2611){\makebox(0,0)[lb]{$u$}}
\put(1551,-3136){\makebox(0,0)[lb]{$-$}}
\end{picture}
\caption{Standard Feedback Configuration} \label{fig1} }
\end{figure}

The transfer function of the plant is
\[
G(s) = \f{q}{s - p}
\]
where $p$ and $q$ are uncertain  parameters. These parameters are assumed as
independent Gaussian
random variables
with density
$\cal{N}(q_0, \si_q)$ and $\cal{N}(p_0, \si_p)$ respectively and
\be
\la{con_a}
p_0 < 0, \qqu  q_0 > 0
\ee
For the worst case design it is assumed that $(q,p) \in \cal{B} (r)$ where
\[
\cal{B} (r) = \{ (x, y) \mid \; |x - q_0| \leq r, \qu |y - p_0| \leq r \}
\]
is the uncertainty bounding set with radius $r > 0$.
We use $P^\mrm{Box}$ to denote $\Pr \{ (q, p) \in \cal{B} (r) \}$, i.e.,
the coverage probability of the uncertainty bounding set.  It can be shown
that
\[
P^\mrm{Box} = \mrm{erf} \li ( \f{r}{\sq{2} \; \si_p} \ri ) \mrm{erf} \li (
\f{r}{\sq{2} \; \si_q} \ri )
\]
where
\[
\mrm{erf}(x)  \DEF \f{2}{\sq{\pi}} \int_0^x e^{-t^2} d t.
\]
Hence, by increasing the radius $r$, one can reduce the uncertainty
modeling error.
\bsk

Consider two controllers
\[
C_A = \f{ K_A } { s + a }, \qqu a > 0
\]
and
\[
C_B = K_B.
\]
We assume further that
\be
\la{con_b}
1 < K_B < \f{K_A}{a}.
\ee

\bsk

When using controller $C_A$, the closed-loop transfer function is given by
\[
T = \f{ \f{ K_A } { s + a } \f{q}{s - p} } { 1 + \f{ K_A } { s + a } \f{q}{s
- p} }
= \f{ K_A  q } { s^2 + (a - p)s  +  K_A  q  - a p}
\]
which is stable if and only if
\[
p < a, \qqu p < \f{K_A}{a} q.
\]
The deterministic stability margin of the system is given by
\[
\ro_A  = \sup \li \{  r > 0 \mid p_0 + r < \f{K_A}{a} (q_0 - r),  \qu p_0 +
r < a \ri \}
\]
which can be simplified as
\[
\ro_A = \min \li ( \f{ K_A \; q_0 - a \; p_0 } {K_A + a}, \; a - p_0 \ri ).
\]
When using controller $C_B$, the closed-loop transfer function is given by
\[
T = \f{ K_B \; q } { s + K_B \; q - p }
\]
which is stable if and only if
\[
p < K_B \; q.
\]
The deterministic stability margin of the system is given by
\bee
\ro_B & = & \sup \li \{ r > 0 \mid p_0 + r < K_B (q_0 - r) \ri \}\\
& = & \f{ K_B \; q_0 - p_0 } {K_B + 1}.
\eee
For uncertainty bounding set with radius $r \in ( \ro_B, \; \ro_A)$, i.e.,
\be
\la{con_c}
\f{ K_B \; q_0 - p_0 } {K_B + 1} < r < \min \li ( \f{ K_A \; q_0 - a \; p_0
} {K_A + a}, \; a - p_0 \ri ),
\ee
controller $B$ may make the system unstable while
controller $A$ robustly stabilizes the system.  More specifically,
controller $B$
can only stabilize a proportion of the family of uncertain plants.
Such a proportion, denoted by $\bb{P}(r)$, is referred to as {\em proportion
of stability},
which is computed as the ratio of the volume of the set
of parameters making the system stable to the total volume of the
uncertainty box, i.e.,
\[
\bb{P}(r) = \f{ \mrm{vol} \{ (q, p) \in \cal{B}(r) \mid
\tx{The system is stable for $(q, p)$} \}  } { \mrm{vol} \{ \cal{B}(r) \} }.
\]
Here ``$\mrm{vol}$'' denotes the Lebesgue measure.
More details are given in Appendix A where we show that the proportion of
instability for controller $B$ is given by
\be
\la{pro}
\bb{P}^B (r) = \bec 1 & \tx{for} \qu 0 < r < \ro_B;\\
1 - \f{ K_B \li ( r + \f{p_0 + r}{K_B} - q_0 \ri )^2  } { 8 r^2} & \tx{for}
\qu \ro_B \leq r \leq \ro_B^*;\\
\f{1}{2} - \f{  \f{p_0}{K_B} - q_0  } { 2 r } & \tx{for} \qu  r > \ro_B^*
\eec
\ee
with
\[
\ro_B^* = \f{ K_B \; q_0 - p_0 } {K_B - 1}.
\]
For an uncertainty bounding set with radius $r \in ( \ro_B, \; \ro_A)$,
controller $B$ is actually a probabilistic
controller because its proportion of stability is strictly less than $1$.
Obviously, controller $A$ is a worst-case controller and is naturally
considered to be more reliable than the probabilistic controller $B$.
However, the following exact computation of probabilities of instability for
both controllers
reveals that the worst-case controller can actually be substantially more
risky than the probabilistic controller.

\bsk

We use $P^{C_A}$ to denote $\Pr \{ \tx{Controller $C_A$ de-stabilizes the
system}  \}$.  We have derived an exact expression as
\be
\la{exact}
P^{C_A}  = \f{1}{2 \pi} \li [ \int_{\se = 0}^{\se^*} \exp \li ( - \f{u^2}{2
\sin^2 \se} \ri ) d \se +
\int_{\se = \se^* - \arctan(k)}^{  \pi  } \exp \li ( -
\f{w^2}{2 \sin^2 \se} \ri ) d \se \ri ]
\ee
where
\[
u = \f{a - p_0}{\si_p} > 0, \qqu  v = \f{ \f{K_A}{a} q_0 - p_0 }{ \si_p} >
0, \qqu k = \f{K_A}{a} \f{\si_q}{\si_p}, \qqu w = \f{v}{ \sq{1 + k^2} }
\]
and
\[
\se^* = \arctan \li(  \f{k u}{u - v} \ri) \in (0, \pi).
\]
For a proof of formula (\ref{exact}), see Appendix B.

\bsk
We use $P^{C_B}$ to denote $\Pr \{ \tx{Controller $C_B$ de-stabilizes the
system}  \}$.  We have derived an exact expression as
\be
\la{exact2}
P^{C_B} =  \f{1}{2} - \f{1}{2} \; \mrm{erf} \li ( \f{K_B  q_0 - p_0} {
\sq{2(\si_p^2 + K_B^2 \si_q^2)}  }   \ri ).
\ee
For a proof of formula (\ref{exact2}), see Appendix C.

A sufficient (but not necessary) condition for the worst-case controller to
be more
risky
than the probabilistic controller (i.e., $P^{C_A} > P^{C_B}$) is
\be
\la{suf}
1 + \li ( \f{ K_B \; \si_q } { \si_p } \ri )^2 < \li( \f{K_B q_0 - p_0  } {
a -
p_0 } \ri )^2,
\ee
which can be easily satisfied.  For a derivation of condition (\ref{suf}),
see
Appendix D.

The boundary of stability is shown in Figure \ref{fig4}.

\begin{figure}[htbp]
\centerline{\psfig{figure=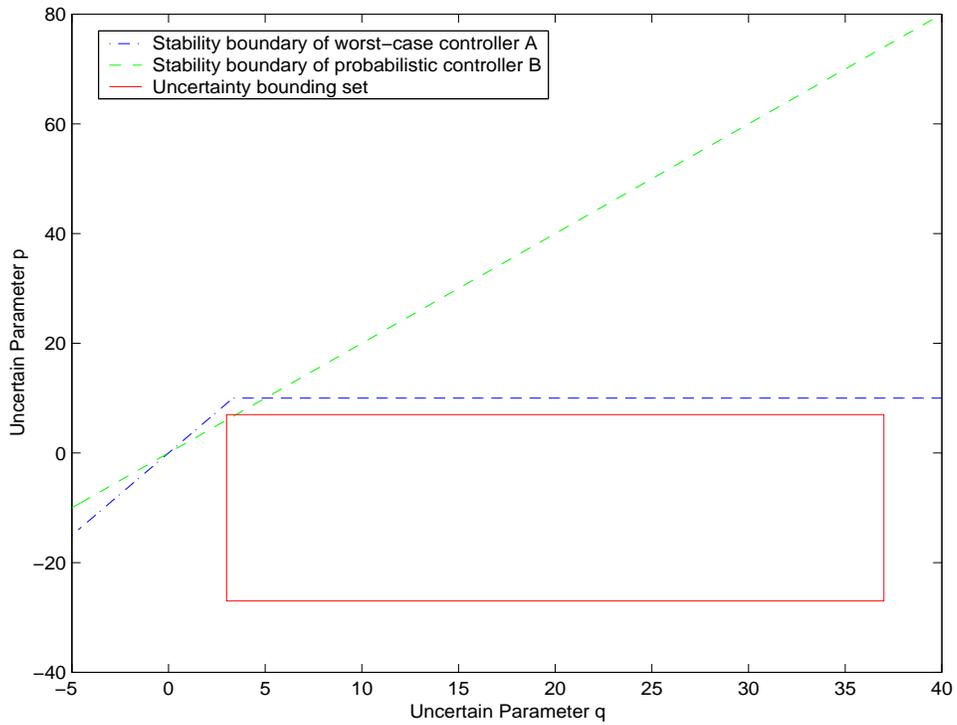, height=3.8in, width=5.0in
}} \caption{Comparison of worst-case controller and probabilistic controller
($r = 17, \; a = 10, \; p_0 = -10,  \; q_0 = 20, \; K_A = 30, \; K_B = 2, \;
\si_p = 10, \; \si_q = 5$)} \la{fig4}
\end{figure}

In Table $1$ we compute the ratio of probabilities of being unstable for
both the worst-case and
probabilistic designs for several situations. The results show that a
worst-case controller can
be thousands of
times more risky than a probabilistic controller.  Granted that this
is a simple first order system but our contention is that  if it happens
even for this simple case then the situation can (easily) happen for highly
complicated systems.

\begin{table}
\caption{Comparison of risks ($\ro_B < r < \ro_A$)} \label{table_example}
\begin{center}
\begin{tabular}{|c||c||c||c||c||c||c||c||c|c|c|c|}
\hline $a$ & $r$ & $p_0$ & $q_0$ & $\si_p$ & $\si_q$ & $K_A$ & $K_B$ &
$P^\mrm{Box}$ &
$\bb{P}^B(r)$ & $P^{C_A}$ & $\f{ P^{C_A} } { P^{C_B} }$\\
\hline $10$ & $17$ & $-10$ & $20$ & $10$ & $5$ & $30$ & $2$ & $0.91$ &
$0.9998$ & $2.3
\times 10^{-2}$ & $112$\\
\hline $40$ & $49$ & $-10$ & $50$ & $20$ & $10$ & $4000$ & $10$ & $0.99$ &
$0.9956$ & $6.2
\times 10^{-3}$ & $2.2 \times 10^{4}$\\
\hline

\end{tabular}
\end{center}
\end{table}
It has been the consensus in the field that guarantees with certainty are
often required for stability and performance in control,
while the probabilistic design is a viable approach when only probabilistic
guarantees are required (see, e.g., lines 16-25, page 4652 of \cite{sc}).
From our general arguments and the concrete example, we can see that, in a
strict sense,
``guarantees with certainty'' are only possible within the modeled uncertainty
bounding set.
Our contention is that such worst-case guarantees may not imply better robustness than
probabilistic
guarantees because of the fact the uncertainty bounding set may not include
all instances of uncertainty.

\sect{Conclusion}

In this paper, we demonstrate that the deterministic worst-case robust
control design
does not necessarily imply a risk free solution and that, in fact,
it can be more risky than a probabilistic controller.  In the final
instance, every
design and analysis is subject to a level of risk.  The goal of design
should be to
make the risk acceptable, instead of assuming that it can make it vanish.

It has been acknowledged that probabilistic methods overcome
the computational complexity and conservatism of worst-case approach at the
expense of a probabilistic risk (\cite{Cal34}, p. 1908). In fact has been remarked by Kargonekar and Tikku that {\it if one is willing to
draw conclusions
with a high degree of confidence, then the computational complexity
decreases dramatically} (\cite{KT}, p. 3470).
More formally, let $\ep, \de \in (0,1)$ and define a system to be {\it
$\ep$-non-robust} if \[
\f{ \mrm{vol} ( \{q \in  \mathcal{B} \mid  \mathbf{P} \; \mrm{is \; violated \; for}
\; q \}) } { \mrm{vol} ( \mathcal{B} ) }\geq \ep.
\]
Then one can detect any $\ep$-non-robust system with probability greater
than
$1 - \de$ based on $N > \f{ \ln \f{1}{\de} } {\ln \f{1}{1-\ep}  }$ i.i.d.
simulations \cite{KT, TD}.

Recently, Barmish et. al. pointed out that {\it if one is
willing to accept some small risk probability of performance violation, it is
often possible to expand the radius of allowable uncertainty by a considerable amount
beyond that provided by the classical robustness theory }(\cite{BLT}, p. 853).

Putting these last two statements in the light of our result make a strong  case for the use of probabilistic
robustness methods.

\bsk

\appendix

\section{Proportion of Instability of Probabilistic Controller}

The set of values uncertainties bounded in $\cal{B} (r)$ which make the
system unstable is
\be
\la{set1}
\cal{B}^\mrm{Bad} (r) = \{ (x, y) \mid \; y \geq K_B \; x, \qu |x - q_0|
\leq r,
\qu |y - p_0| \leq r \}.
\ee
Clearly,
\[
\bb{P}^B(r) = 1 - \f{ \mrm{vol} \{ \cal{B}^\mrm{Bad} (r) \} } { \mrm{vol} \{
\cal{B} (r) \}  }.
\]
We can now compute the proportion of stability $\bb{P}(r)$ for controller
$C_B$ as follows.

\bsk

It can be shown that $\cal{B}^\mrm{Bad} (r) = \emptyset$ for $r < \ro_B$.
Hence
\[
\bb{P}^B(r) = 1, \qqu 0 \leq r < \ro_B.
\]

\bsk

For $\ro_B \leq r \leq \ro_B^*$, we have
\be
\la{set2}
\cal{B}^\mrm{Bad} (r) = \li \{ (x, y) \mid \; p_0 + r \geq y \geq K_B \; x,
\qu
q_0 - r \leq x \ri \}.
\ee
We now prove equation (\ref{set2}).  For notation simplicity, let the set in
the right-hand side of (\ref{set1}) be denoted as $\Se$.  Let the set in
the right-hand side of (\ref{set2}) be denoted as $\Phi$.  Clearly, $\Se
\subseteq \Phi$.
It suffices to show $\Se \supseteq \Phi$.  Let $(x, y) \in \Phi$.
Then $y \geq K_B \; x \geq K_B \; (q_0 - r)$.
It can be verified that $K_B \; (q_0 - r) \geq p_0 - r$ if and only if $r
\leq \ro_B^*$.
Hence $p_0 + r \geq y \geq K_B \; x \geq p_0 - r$.
From inequalities (\ref{con_a}), (\ref{con_b}) and $p_0 + r \geq K_B \; x$,
we obtain $x \leq \f{p_0 + r}{K_B} \leq q_0 + r$.  Therefore, $(x, y) \in
\Se$.
It follows that $\Se \supseteq \Phi$ and thus $\Se = \Phi$.
Observing that $\cal{B}^\mrm{Bad} (r)$ is a triangular domain, we find by a
geometric method
\[
\mrm{vol} \{ \cal{B}^\mrm{Bad} (r) \} = \f{1}{2} K_B \li ( r + \f{p_0 +
r}{K_B} - q_0 \ri )^2.
\]
It follows that
\[
\bb{P}^B(r) = 1 - \f{ K_B \li ( r + \f{p_0 + r}{K_B} - q_0 \ri )^2  } { 8
r^2}, \qqu  \ro_B \leq r \leq \ro_B^*.
\]

\bsk

For $r > \ro_B^*$, we have
\be
\la{set3}
\cal{B}^\mrm{Bad} (r) = \Xi \cup \Psi
\ee
where
\[
\Xi = \li \{ (x, y) \mid \; |y - p_0| \leq r, \qu q_0 - r \leq x \leq \f{p_0
-
r}{K_B} \ri \}
\]
and
\[
\Psi = \li \{ (x, y) \mid \; p_0 + r \geq y \geq K_B \; x, \qu \f{p_0 -
r}{K_B}
< x \ri \}.
\]
We now show $(\ref{set3})$, i.e., $\Se = \Xi \cup \Psi$.
Obviously, $\Se \subseteq \Xi \cup \Psi$.  It suffices to show $\Se
\supseteq \Xi$ and
$\Se \supseteq \Psi$.  Let $(x,y) \in \Xi$.
Then, $K_B \; x \leq K_B \; \li (  \f{p_0 - r}{K_B} \ri ) = p_0 - r \leq y$.
From inequalities (\ref{con_a}), (\ref{con_b}) and $x \leq \f{p_0 -
r}{K_B}$,
we have $x \leq \f{p_0 - r}{K_B} \leq q_0 + r$.  This proves $(x, y) \in
\Se$.
Hence $\Xi \subseteq \Se$.  Now let $(x, y) \in \Psi$.
Then $y \geq K_B \; x > K_B \; \li ( \f{p_0 - r}{K_B} \ri )= p_0 - r$.
It can be shown that $\f{p_0 - r}{K_B} > q_0 - r$ if and only if $r >
\ro_B^*$.  Hence, $x > q_0 - r$.  From inequalities (\ref{con_a}),
(\ref{con_b})
and $K_B \; x \leq p_0 + r$,
we have $x \leq \f{p_0 + r}{K_B} \leq q_0 + r$.
This proves $(x, y) \in \Se$ and thus $\Psi \subseteq \Se$.
Therefore, the proof for $\Se = \Xi \cup \Psi$ is completed.

\bsk

Observing that $\Xi$ is a rectangular domain and $\Psi$ is a triangular
domain, we obtain by a geometric argument
\[
\mrm{vol} \{ \cal{B}^\mrm{Bad} (r) \} = 2 r \li ( r + \f{p_0}{K_B} - q_0 \ri
)
\]
and thus
\[
\bb{P}^B (r) = 1 - \f{ 2 r \li ( r + \f{p_0}{K_B} - q_0 \ri ) } { 4 r^2 } =
\f{1}{2} - \f{  \f{p_0}{K_B} - q_0  } { 2 r }, \qqu  r > \ro_B^*.
\]

\bsk

\section{Probability of Instability of The Worst-case Controller }

Now we derive an exact expression for $\Pr \{ \tx{Controller $C_A$
stabilizes the system}  \}$.  Note that
\bee
&   & \Pr \{ \tx{Controller $C_A$ stabilizes the system}  \}\\
& = & \Pr \li \{p < a, \qu p < \f{K_A}{a} q  \ri \}\\
& = & \f{1}{2 \pi \si_p \si_q} \; \int_{ \li \{p < a, \qu p < \f{K_A}{a} q
\ri \}  } \; \exp \li( -
\f{(p - p_0)^2 } { 2
\si_p^2 } \ri) \; \exp \li( -
\f{(q - q_0)^2 } { 2
\si_q^2 } \ri) d p \; dq \\
\eee
Introducing new variables
\[
y = \f{p - p_0  } { \si_p }, \qqu x = \f{q - q_0  } { \si_q },
\]
we have
\bee
&   & \Pr \{ \tx{Controller $C_A$ stabilizes the system}  \}\\
& = & \f{1}{ 2 \pi} \; \int_{D_{xy}} \; \exp \li ( - \f{x^2}{2} \ri ) \;
\exp \li ( - \f{y^2}{2} \ri ) d x \; d y
\eee
where
\bee
D_{xy} & = & \li \{ (x, y) \mid \si_p y + p_0 < a, \qu   \si_p y + p_0 <
\f{K_A}{a} (\si_q x + q_0) \ri \}\\
& = & \li \{ (x, y) \mid y < u, \qu y < k x + v \ri \}.
\eee
Therefore, it suffices to compute
\[
I = \int_{D_{xy}} \; \exp \li ( - \f{x^2 + y^2}{2} \ri ) \;  d x \; d y.
\]
Introducing polar coordinates
\be
\la{polar}
x = \ro \cos \se, \qqu y = \ro \sin \se,
\ee
we have
\[
I = \int_{D_{\ro \se}} \; \exp \li ( - \f{\ro^2}{2} \ri ) \; \ro \; d \ro
\; d \se
\]
where
\[
D_{\ro \se} = \li \{ (\ro, \se) \mid \ro \sin \se  < u, \qu \ro \sin \se < k
\; \ro \cos \se + v,
\qu \ro \in [0, \iy), \qu \se \in [0, 2 \pi) \ri \}.
\]
We compute the integration
\[
\wh{I} = \int_{\wh{D}_{\ro \se}} \; \exp \li ( - \f{\ro^2}{2} \ri ) \; \ro
\; d \ro  \; d \se
\]
over the complement set
\[
\wh{D}_{\ro \se} = \li \{ (\ro, \se) \mid \ro \sin \se  \geq u \qu \tx{or}
\qu \ro \sin \se \geq k \; \ro \cos \se + v,
\qu \ro \in [0, \iy), \qu \se \in [0, 2 \pi) \ri \}.
\]
We claim that $\wh{D}_{\ro \se}$ can be partitioned as two subsets
\[
A_{\ro \se} = \li \{ (\ro, \se) \mid \ro \sin \se  \geq u, \qu \ro \in [0,
\iy), \qu \se \in [0, \se^*) \ri \}
\]
and
\[
B_{\ro \se} = \li \{ (\ro, \se) \mid  \ro \sin \se \geq k \; \ro \cos \se +
v,
\qu \ro \in [0, \iy), \qu \se \in [\se^*, 2 \pi) \ri \}
\]
such that
\[
\wh{D}_{\ro \se} = A_{\ro \se} \cup B_{\ro \se}, \qqu A_{\ro \se} \cap
B_{\ro \se} = \emptyset.
\]

\bsk

In the proof of the claim, we shall recall that $(x, y)$ is related to
$(\ro, \se)$ by the polar transform (\ref{polar}).
To show $\wh{D}_{\ro \se} = A_{\ro \se} \cup B_{\ro \se}$, we need to
consider
three cases: Case (a): $u > v$; Case (b): $u < v$; Case (c): $u = v$.

\bsk

Let $k^* = \tan \se^* = \f{k u}{u - v}$.  In Case (a), since $\se < \se^*$
if and only if $y < k^* x$,
it suffices to show the following two statements:

\bed
\item [(a-1)] $y > u$ if $y < k^* x, \qu y \geq k x + v$;

\item [(a-2)] $y \geq k x + v$ if $y \geq u, \qu y \geq k^* x$.

\eed

To show statement (a-1), observing that,  as a direct consequence of $y
< k^* x$ and $y \geq k x + v$, we have $k^* x > k x + v$, leading to $x >
\f{v}{k^* - k}$.
Therefore, $y \geq k x + v > k \f{v}{k^* - k} + v = u$.
To show statement (a-2), we proceed by a contradiction method.
Suppose $y < k x + v$.  Then $k^* x < k x + v$, which implies $x < \f{v}{k^*
- k}$.
On the other hand, we have $k x + v > u$, leading to $x > \f{u- v}{k}$.  It
follows that
$\f{u- v}{k} < \f{v}{k^* - k} = \f{u- v}{k}$, which is a contradiction.

\bsk

In Case (b), since $\se < \se^*$ if and only if $y > k^* x$, it suffices to
show the following two statements:

\bed
\item [(b-1)] $y > u$ if $y > k^* x, \qu y \geq k x + v$;

\item [(b-2)] $y \geq k x + v$ if $y \geq u, \qu y \leq k^* x$.

\eed

To show statement (b-1),  observing that $k^* < 0$ because $u < v$.
Hence $x > \f{y}{k^*}$.
On the other hand $x \leq \f{y - v}{k}$.  Hence, $\f{y}{k^*} < \f{y -
v}{k}$, leading to $ y > u$.  To show statement (b-2), we proceed by a
contradiction method.
Suppose $y < k x + v$.  Then $u < k x + v$, i.e., $x > \f{u-v}{k}$.
Consequently, $y \leq k^* x < k^* \f{u-v}{k} = u$, which is a contradiction.

\bsk

In Case (c), since $\se < \se^* = \f{\pi}{2}$ if and only if $x > 0$, it
suffices to show the following two statements:

\bed

\item [(c-1)] $y > u$ if $x > 0, \qu y \geq k x + v$;
\item [(c-2)] $y \geq k x + v$ if $y \geq u, \qu x \leq 0$.

\eed

Statement (c-1) can be shown by observing that, if $x > 0, \; y \geq k x
+ v$, then $y > v = u$.  To show statement (c-2), we proceed by a
contradiction method.
Suppose $y < k x + v$.  Then $u < k x + v$, i.e., $x > \f{u-v}{k} = 0$,
which is a contradiction.

\bsk

It can be shown that
\be
\la{sim1}
A_{\ro \se} = \li \{ (\ro, \se) \mid \ro  \geq \f{u}{ \sin \se },  \qu \se
\in (0, \se^*) \ri \}
\ee
and $B_{\ro \se} = \li \{ (\ro, \se) \mid  \ro \sin \se \geq k \; \ro \cos
\se +
v, \qu \ro \in [0, \iy), \qu \se^* \leq \se < \pi + \arctan(k)  \ri \}$.
Moreover, we can further simplify $B_{\ro \se}$ as
\be
\la{sim}
B_{\ro \se} = \li \{ (\ro, \se) \mid \ro \geq \f{v} { \sin \se - k \cos \se
},
\qu \se^* \leq \se < \pi + \arctan(k)  \ri \}.
\ee
To show (\ref{sim}), it suffices to show
\be
\la{ineq1}
\sin \se - k \cos \se \geq 0 \qu \tx{for} \qu \se^* \leq \se \leq \f{\pi}{2}
\ee
and
\be
\la{ineq2}
\sin \se - k \cos \se \geq 0 \qu \tx{for} \qu  \f{\pi}{2} < \se < \pi +
\arctan(k).
\ee
To show (\ref{ineq1}), one needs to observe that $ \se^* \leq \se <
\f{\pi}{2}$ implies
\[
\cos \se > 0, \qu u \geq v, \qu \tan \se \geq \tan (\se^*) = \f{k u  } {u -
v } \geq k
\]
and consequently,  $\sin \se - k \cos \se = \cos \se \; (\tan \se - k ) \geq
0$.
To show (\ref{ineq2}), one needs to observe that $\f{\pi}{2} < \se < \pi +
\arctan(k)$ implies
\[
\cos \se < 0, \qu \tan \se < \tan (\pi + \arctan(k)) = k
\]
and consequently,  $\sin \se - k \cos \se = \cos \se \; (\tan \se - k ) >
0$.   By (\ref{sim1}) and (\ref{sim}), we have
\bee
\wh{I} & = & \int_{A_{\ro \se}} \; \exp \li ( - \f{\ro^2}{2} \ri ) \; \ro \;
d \ro  \; d \se +
\int_{B_{\ro \se}} \; \exp \li ( - \f{\ro^2}{2} \ri ) \; \ro \; d \ro  \; d
\se\\
& = & \int_{\se = 0}^{\se^*} \; \int_{\ro = \f{u}{ \sin \se }}^\iy \exp \li
( - \f{\ro^2}{2} \ri ) \; \ro \; d \ro  \; d \se\\
&   &  + \int_{\se = \se^*}^{ \pi + \arctan(k) }  \int_{\ro =  \f{v} { \sin
\se - k \cos \se } }^\iy\;
\exp \li ( - \f{\ro^2}{2} \ri ) \; \ro \; d \ro  \; d \se\\
& = & \int_{\se = 0}^{\se^*} \exp \li ( - \f{u^2}{2 \sin^2 \se} \ri ) d \se
+ \int_{\se = \se^*}^{ \pi + \arctan(k) } \exp \li ( - \f{v^2}{2 (\sin \se -
k \cos \se)^2} \ri ) d \se\\
& = & \int_{\se = 0}^{\se^*} \exp \li ( - \f{u^2}{2 \sin^2 \se} \ri ) d \se
+ \int_{\se = \se^* - \arctan(k)}^{ \pi } \exp \li ( -
\f{w^2}{2 \sin^2 \se} \ri ) d \se.
\eee
Finally,
\[
\Pr \{ \tx{Controller $C_A$ stabilizes the system}  \} = 1 - \f{1}{2 \pi}
\wh{I}.
\]
This completes the proof of formula (\ref{exact}).

\bsk

\section{Probability of Instability of Probabilistic Controller}

Define
\[
S_{pq} = \li \{ (p,q) \mid p \in \bb{R}, \qu q \in \bb{R}, \qu p < K_B \; q
\ri \}.
\]
Then \bee
&   & \Pr \{ \tx{Controller $C_B$ stabilizes the system}  \}\\
& = & \Pr \{ (p,q) \in  S_{pq}\}\\
& = & \f{1}{2 \pi \si_p \si_q} \int_{S_{pq}} \exp \li( -
\f{(p - p_0)^2 } { 2 \si_p^2 } \ri) \exp \li( - \f{(q - q_0)^2 } { 2
\si_q^2 } \ri) d p d q. \eee Note that there exists a one-to-one
mapping between $S_{pq}$ and
\[
S_{xy} = \{ (x, y) \mid x \in \bb{R}, \; y < 0  \}
\]
so that
\[
p = K_B \; x + y,  \qqu q = x.
\]
Therefore,
\bee
&   & \Pr \{ \tx{Controller $C_B$ stabilizes the system}  \}\\
& = &  \f{1}{2 \pi \si_p \si_q} \int_{S_{xy}} \exp \li( -
\f{(K_B x + y - p_0)^2 } { 2 \si_p^2 } \ri)
\exp \li( - \f{(x - q_0)^2 } { 2 \si_q^2 } \ri) d x d y\\
& = & \f{1}{2 \pi \si_p \si_q} \int_{S_{xy}} \exp \li(  - \f{
\si_p^2 + K_B^2 \si_q^2  } { 2 \si_p^2 \si_q^2 } x^2  + \f{q_0 \si_p^2 - K_B
\si_q^2 (y - p_0) }
{\si_p^2 \si_q^2} x  - \f{(y - p_0)^2 } { 2
\si_p^2 } - \f{q_0^2} { 2 \si_q^2 }  \ri) d x d y\\
& = & \f{1}{2 \pi \si_p \si_q} \int_{-\iy}^0 \exp \li ( \f{ [q_0 \si_p^2 -
K_B \si_q^2 (y - p_0) ]^2 }
{ 2 \si_p^2 \si_q^2 (\si_p^2 + K_B^2 \si_q^2)}  - \f{(y - p_0)^2 } { 2
\si_p^2 } - \f{q_0^2} { 2 \si_q^2 }  \ri) \\
&   & \times \int_{-\iy}^\iy
\exp \li(  - \f{\si_p^2 + K_B^2 \si_q^2} { 2 \si_p^2 \si_q^2 }
\li [ x  - \f{q_0 \si_p^2 - K_B \si_q^2 (y - p_0) } {\si_p^2 + K_B^2
\si_q^2} \ri ] ^2 \ri ) d x \; d y.
\eee
Using the fact
\[
\int_{-\iy}^\iy e^{-\al (x - \mu)^2} dx = \sq{ \f{\pi}{\al}  }  \qqu \fa
\al > 0, \qu \fa \mu \in (- \iy, \iy),
\]
we have \bee
&   & \Pr \{ \tx{Controller $C_B$ stabilizes the system}  \}\\
& = & \f{1}{2 \pi \si_p \si_q} \int_{-\iy}^0 \; \exp \li ( \f{ [q_0 \si_p^2
- K_B \si_q^2 (y - p_0) ]^2 }
{ 2 \si_p^2 \si_q^2 (\si_p^2 + K_B^2 \si_q^2)}  - \f{(y - p_0)^2 } { 2
\si_p^2 } - \f{q_0^2} { 2 \si_q^2 }  \ri) d y  \;
    \times \sq{ \f{ 2 \pi \si_p^2 \si_q^2 } { \si_p^2 + K_B^2 \si_q^2 } }\\
& = & \f{1}{2 \pi } \; \sq{ \f{ 2 \pi  } { \si_p^2 + K_B^2 \si_q^2 } } \;
\int_{-\iy}^{-p_0} \exp \li ( \f{ [q_0 \si_p^2 - K_B \si_q^2 z ]^2 }
{ 2 \si_p^2 \si_q^2 (\si_p^2 + K_B^2 \si_q^2)}  - \f{z^2 } { 2 \si_p^2 } -
\f{q_0^2} { 2 \si_q^2 }  \ri) d z.
\eee
It can be verified that
\[
\f{ [q_0 \si_p^2 - K_B \si_q^2 z ]^2 }
{ 2 \si_p^2 \si_q^2 (\si_p^2 + K_B^2 \si_q^2)}  - \f{z^2 } { 2 \si_p^2 } -
\f{q_0^2} { 2 \si_q^2 } =
- \f{ \li( z+ K_B \; q_0 \ri)^2 } { 2(\si_p^2 + K_B^2 \si_q^2) }.
\]
Hence
\bee
&   & \Pr \{ \tx{Controller $C_B$ stabilizes the system}  \}\\
& = & \f{1}{2 \pi } \; \sq{ \f{ 2 \pi  } { \si_p^2 + K_B^2 \si_q^2 } } \;
\int_{-\iy}^{-p_0} \exp \li (  - \f{ \li( z+ K_B \; q_0 \ri)^2 } { 2(\si_p^2
+ K_B^2 \si_q^2) }  \ri ) dz\\
& = & \f{1}{\sq{2 \pi}} \int_{-\iy}^{ \f{K_B  q_0 - p_0} {  \sq{\si_p^2 +
K_B^2 \si_q^2}  }  } e^{-\f{z^2}{2}} dz\\
& = &  \f{1}{2} + \f{1}{2} \; \mrm{erf} \li ( \f{K_B  q_0 - p_0} {
\sq{2(\si_p^2 + K_B^2 \si_q^2)}  }   \ri ). \eee
It follows that
\bee
P^{C_B} & = & 1 - \Pr \{ \tx{Controller $C_B$ stabilizes the system}  \}\\
& = & \f{1}{2} - \f{1}{2} \; \mrm{erf} \li ( \f{K_B  q_0 - p_0} {
\sq{2(\si_p^2 + K_B^2 \si_q^2)}  }   \ri ).
\eee
This completes the proof of formula (\ref{exact2}).

\bsk

\section{Derivation of A Sufficient Condition}

Note that
\bee
P^{C_A} & = & 1 - \Pr \li \{p < a, \qu p < \f{K_A}{a} q  \ri \}\\
& \geq & 1 - \Pr \{ p < a \}\\
& = & 1 - \f{1}{\sq{2 \pi} \si_p} \int_{- \iy}^a \exp \li( -
\f{(p - p_0)^2 } { 2
\si_p^2 } \ri) dp\\
& = & 1 - \f{1}{\sq{2 \pi} } \int_{- \iy}^{\f{a - p_0}{\si_p}} \exp \li( -
\f{x^2 } { 2  } \ri) dx\\
& = & \f{1}{2} - \f{1}{2} \; \mrm{erf} \li ( \f{a - p_0}{ \sq{2}
\si_p} \ri ). \eee
Therefore, for $P^{C_A} > P^{C_B}$, it suffices to have
\[
\f{1}{2} - \f{1}{2} \; \mrm{erf} \li ( \f{a - p_0}{ \sq{2}
\si_p} \ri ) > \f{1}{2} - \f{1}{2} \; \mrm{erf} \li ( \f{K_B  q_0 - p_0} {
\sq{2(\si_p^2 + K_B^2 \si_q^2)}  }   \ri ),
\]
i.e.,
\[
\mrm{erf} \li ( \f{a - p_0}{ \sq{2}
\si_p} \ri ) < \mrm{erf} \li ( \f{K_B  q_0 - p_0} {  \sq{2(\si_p^2 + K_B^2
\si_q^2)}  }   \ri ).
\]
Since $\mrm{erf}(.)$ is a monotone increasing function, we have
\[
\f{a - p_0}{ \sq{2} \si_p}  < \f{K_B  q_0 - p_0} {  \sq{2(\si_p^2 + K_B^2
\si_q^2)}  },
\]
which is equivalent to
\[
1 + \li ( \f{ K_B \; \si_q } { \si_p } \ri )^2 < \li ( \f{K_B q_0 - p_0  } {
a -
p_0 } \ri )^2.
\]
This completes the proof of condition (\ref{suf}).

\bsk


\begin{thebibliography}{10}

\bibitem{bai}E. W. BAI, R. TEMPO,  AND M. FU,
``Worst-case properties of the uniform distribution and randomized
algorithms for robustness analysis,'' {\it Mathematics of Control,
Signals and Systems}, vol. 11, pp.183-196, 1998.

\bibitem{BLT} B. R. BARMISH, C. M. LAGOA, AND R. TEMPO,
``Radially truncated uniform distributions for probabilistic
robustness of control systems,'' {\it Proc. of American Control
Conference}, pp. 853-857, Albuquerque, New Mexico, June 1997.


\bibitem{BL}
B. R. BARMISH AND C. M. LAGOA, ``The uniform distribution: a
rigorous justification for its use in robustness analysis,'' {\it
Mathematics of Control, Signals and Systems}, vol. 10, pp.
203-222, 1997.

\bibitem{Barmish1} B. R. BARMISH AND P. S. SHCHERBAKOV, ``A dilation method
for robustness problems with nonlinear parameter dependence,''
{\it Proc. of American Control Conference}, pp. 3834-3839, Denver,
2003.

\bibitem{Barmish2} B. R. BARMISH AND P. S. SHCHERBAKOV, ``On avoiding
vertexization of robustness problems: The approximate feasibility
concept,'' {\it IEEE Transactions on Automatic Control}, vol. 42,
pp. 819-824, 2002.


\bibitem{CDT1}
G. CALAFIORE, F. DABBENE, AND R. TEMPO, ``Randomized algorithms
for probabilistic robustness with real and complex structured
uncertainty,'' {\em IEEE Transaction on Automatic Control}, vol.
45, pp. 2218-2235, 2000.

\bib{Cal2} G. CALAFIORE AND M. C. CAMPI, ``Uncertain convex programs:
randomized solutions and confidence levels,'' to appear in {\it
Mathematical Programming}, 2004.

\bib{Cal33} G. CALAFIORE AND F. DABBENE, ``A probabilistic framework for
problems with real structured uncertainty in systems and
control,'' {\it Automatica}, vol. 38, pp. 1265-1276, 2002.

\bib{Cal3} G. CALAFIORE AND B. T. POLYAK, ``Fast algorithms for exact and
approximate feasibility of robust LMIs,'' {\it IEEE Transaction on
Automatic Control}, vol. 46, pp. 1755-1759, 2001.

\bib{Cal34} G. CALAFIORE AND F. DABBENE, AND R. TEMPO, ``Randomized
algorithms in robust control,''
{\it Proceedings IEEE Conference on Decision and Control}, pp. 1908-1913,
Maui, December 2003.

\bibitem{C0} X. CHEN, K. ZHOU, AND J. ARAVENA, ``Fast construction of
robustness degradation function,'' {\it SIAM Journal on Control
and Optimization}, vol. 42, pp. 1960-1971, 2004.

\bibitem{C1} X. CHEN, K. ZHOU, AND J. ARAVENA, ``Fast universal
algorithms for robustness analysis,'' {\it Proceedings IEEE
Conference on Decision and Control}, pp. 1926-1931, Maui, December
2003.

\bibitem{C2} X. CHEN AND K. ZHOU,
``Constrained robustness analysis and synthesis by randomized
algorithms,'' {\em IEEE Transaction on Automatic Control}, vol.
45, pp. 1180-1186, 2000.

\bib{Dan}  D. P. FARIAS AND B. V. ROY, ``On constraint sampling in the
linear programming approach to approximate linear programming,''
{\it Proceedings IEEE Conference on Decision and Control}, pp.
2441-2446, Maui, December 2003.

\bib{FU} Y. FUJISAKI, F. DABBENE, AND R. TEMPO, ``Probabilistic robust
design of LPV control systems,'' {\it Automatica}, vol. 39, pp.
1323-1337, 2003.

\bib{FU1} Y. FUJISAKI AND Y. KOZAWA, ``Probabilistic Rrobust controller
design: probable near minimax
value and randomized algorithms,'' {\it Proceedings IEEE
Conference on Decision and Control}, pp. 1938-1943, Maui, December
2003.

\bib{PT} P. F. HOKAYEM AND C. T. ABDALLAH, ``Quasi-Monte Carlo methods in
robust
control design,'' {\it Proceedings IEEE Conference on Decision and
Control}, pp. 2435-2440, Maui, December 2003.

\bib{Kan} S. KANEV, B. De SCHUTTER, AND M. VERHAEGEN, ``An ellipsoid
algorithm for probabilistic robust controller design,'' {\it
Systems and Control Letters}, vol. 49, pp. 365-375, 2003.


\bib{Kan2} S. KANEV AND M. VERHAEGEN,
``Robust output-feedback integral MPC: A probabilistic approach,''
{\it Proceedings IEEE Conference on Decision and Control}, pp.
1914-1919, Maui, December 2003.

\bibitem{KT}
P. P. KHARGONEKAR AND A. TIKKU, ``Randomized algorithms for robust
control analysis and synthesis have polynomial complexity,'' {\it
Proceedings of Conference on Decision and Control}, pp. 3470-3475,
Kobe, Japan, December 1996.

\bib{Ko} V. KOLTCHINSKII, C.T. ABDALLAH, M. ARIOLA, P. DORATO, AND D.
PANCHENKO, ``Improved sample complexity estimates for statistical
learning control of uncertain systems.,'' {\it IEEE Transactions
on Automatic Control}, vol. 46, pp. 2383-2388, 2000.

\bibitem{BLT2} C. M. LAGOA, ``Probabilistic enhancement  of classic
robustness margins: A class of none symmetric distributions,''
{\it Proc. of American Control Conference}, pp. 3802-3806,
Chicago, Illinois, June 2000.

\bib{La2} C. M. LAGOA, X. LI, AND M. SZNAIER, ``On the
design of robust controllers for arbitrary uncertainty
structures,'' to appear in {\it IEEE Transactions on Automatic
Control}.

\bib{La3} C. M. LAGOA, X. LI,  M. C.
MAZZARO, AND M. SZNAIER, ``Sampling random transfer functions,''
{\it Proceedings IEEE Conference on Decision and Control}, pp.
2429-2434, Maui, December 2003.

\bib{La4} C. M. LAGOA AND B. R. BARMISH, ``Distributionally robust Monte
Carlo simulation: A tutorial survey,'' {\it Proceedings of the
IFAC World Congress}, pp. 1327-1338, 2002.

\bibitem{MS}
C. MARRISON AND R. F. STENGEL, ``Robust control system design
using random search and genetic algorithms,'' {\em IEEE
Transaction on Automatic Control}, vol. 42, pp. 835-839, 1997.

\bib{nem} A. NEMIROVSKII, ``On tractable approximations of randomly
perturbed convex constraints,'' {\it Proceedings IEEE Conference
on Decision and Control}, pp. 2419-2422, Maui, December 2003.

\bib{Oi} Y. OISHI AND H. KIMURA, ``Randomized algorithms to solve
parameter-dependent linear matrix inequalities and their
computational complexity,'' {\it Proceedings IEEE Conference on
Decision and Control}, pp. 2025-2030, 2001.

\bib{Oi2} Y. OISHI, ``Probabilistic design of a robust state-feedback
controller based
on parameter-dependent Lyapunov functions,'' {\it Proceedings IEEE
Conference on Decision and Control}, pp. 1920-1925, Maui, December
2003.

\bib{Pol} B. T. POLYAK AND R. TEMPO, ``Probabilistic robust design with
linear quadratic regulators,'' {\it Systems and Control Letters},
vol. 43, pp. 343-353, 2001.

\bibitem{Polyak} B. T. POLYAK AND P. S. SHCHERBAKOV,
``Random spherical uncertainty in estimation and robustness,''
{\em IEEE Transaction on Automatic Control}, vol. 45, pp.
2145-2150, 2000.

\bibitem{RS}
L. R. RAY AND R. F. STENGEL, ``A monte carlo approach to the
analysis of control systems robustness,'' {\it Automatica}, vol.
3, pp. 229-236, 1993.

\bibitem{SB}
S. R. ROSS AND B. R. BARMISH, ``Distributionally robust gain
analysis for systems containing complexity,'' {\it Proceedings of
Conference on Decision and Control}, pp. 5020-5025, Orlando,
Florida, December 2001.

\bib{sc} C. W. SCHERER, ``Higher-order relaxations for robust LMI
problems with verifications for exactness,'' {\it Proceedings IEEE
Conference on Decision and Control}, pp. 4652-4657, Maui, December
2003.

\bib{kov} P. S. SHCHERBAKOV AND B. R. BARMISH, ``On the conditioning of
robustness problems,'' {\it Proceedings IEEE Conference on
Decision and Control}, pp. 1932-1937, Maui, December 2003.

\bibitem{SR}
R. F. STENGEL AND L. R. RAY, ``Stochastic robustness of linear
time-invariant systems,'' {\it IEEE Transaction on Automatic
Control}, vol. 36, pp. 82-87, 1991.

\bib{Tem} R. TEMPO, G. CALAFIORE, AND F. DABBENE, {\it Randomized Algorithms
for Analysis and Control of Uncertain Systems}, Springer-Verlag,
New York, 2004.

\bibitem{TD}
R. TEMPO, E. W. BAI, AND F. DABBENE, ``Probabilistic robustness
analysis: explicit bounds for the minimum number of samples,''
{\em Systems and Control Letters}, vol. 30, pp. 237-242, 1997.


\bibitem {vind} M. VIDYASAGAR AND V. D. BLONDEL,
``Probabilistic solutions to NP-hard matrix problems,'' {\it
Automatica}, vol. 37, pp. 1597-1405, 2001.

\bibitem {vind2} M. VIDYASAGAR, ``Randomized algorithms for robust
controller synthesis using statistical learning theory,'' {\it
Automatica}, vol. 37, pp. 1515-1528, 2001.

\bibitem {vind3} M. VIDYASAGAR, ``Statistical learning theory and randomized
algorithms for control,'' {\it IEEE Control Systems Magazine},
vol. 18, pp. 69-85, 1998.

\bib{Wang} Q. WANG AND R. F. STENGEL, ``Robust control of nonlinear
systems with parametric uncertainty,'' {\it Automatica}, vol. 38,
pp. 1591-1599, 2002.

\bib{Wang2} Q. WANG, ``Probabilistic robust control design of polynomial
vector fields,''
{\it Proceedings IEEE Conference on Decision and Control},  pp.
2447-2452, Maui, December 2003.

\end{thebibliography}
\end{document}